\documentclass{amsart}
\usepackage{amssymb,latexsym}

\newtheorem{theorem}{Theorem}[section]

\newtheorem{proposition}[theorem]{Proposition}
\newtheorem{remark}[theorem]{Remark}
\newtheorem{corollary}[theorem]{Corollary}

\numberwithin{equation}{section}

\begin{document}

\title{On the Cycle Structure of Repeated Exponentiation Modulo a Prime Power}

\author{Min Sha}
\address{Institut de Mathematiques de Bordeaux, Universite Bordeaux 1
, 33405 Talence Cedex, France}
\email{shamin2010@gmail.com}






\begin{abstract}
We obtain some results about the repeated exponentiation modulo
a prime power from the viewpoint of arithmetic dynamical systems.
Especially, we extend two asymptotic formulas about periodic points and tails in the case of modulo a prime to
the case of modulo a prime power.
\end{abstract}

\maketitle



\section{Introduction}                                                                                
For a positive integer $M$, denote by $\mathbb{Z}/M\mathbb{Z}$ the residue ring
of $\mathbb{Z}$ modulo $M$ and $(\mathbb{Z}/M\mathbb{Z})^{*}$ the unit group.
For an integer $k\ge 2$, we consider the following endomorphism of $(\mathbb{Z}/M\mathbb{Z})^{*}$,
\begin{equation}
f: (\mathbb{Z}/M\mathbb{Z})^{*} \to (\mathbb{Z}/M\mathbb{Z})^{*}, x\to x^{k}.
\notag
\end{equation}

For any initial value $x\in (\mathbb{Z}/M\mathbb{Z})^{*}$, we repeat the action of $f$, then we get
a sequence
\begin{equation}
x_{0}=x, x_{n}=x_{n-1}^{k}, n=1,2,3,\cdots.
\notag
\end{equation}
This sequence is known as the power generator of pseudorandom numbers.
Studying such sequences in the cases that $M$ is a prime or a product of two distinct primes,
is of independent interest and is also important for
several cryptographic applications, see \cite{BBS,L}. From the viewpoint of cryptography,
there are numerous results about these sequences, see the papers mentioned in \cite{CS}, more
recently see \cite{E} and its references.

If we view $(\mathbb{Z}/M\mathbb{Z})^{*}$ as a vertex set and draw
a directed edge from $a$ to $b$ if $f(a) = b$, then we get a digraph. There are also many results in
this direction, see \cite{W} and the papers mentioned there, more recently see
\cite{SK1,SK2,SK3,VS}.

As \cite{CS}, in this article we will study $(\mathbb{Z}/M\mathbb{Z})^{*}$ under the action of $f$ from the viewpoint of
arithmetic dynamical systems, where $M$ is a prime power.
Especially we will extend two asymptotic formulas in \cite{CS} to the case of modulo a prime power.

It is easy to see that for any initial value $x\in (\mathbb{Z}/M\mathbb{Z})^{*}$ the corresponding
sequence becomes eventually periodic, that is, for some positive integer $s_{k,M}(x)$
and \textit{tail} $t_{k,M}(x)<s_{k,M}(x)$, the elements $x_{0}=x, x_{1},\cdots, x_{s_{k,M}(x)-1}$ are pairwise distinct
and $x_{s_{k,M}(x)}=x_{t_{k,M}(x)}$. So we can define a \textit{tail function} $t_{k,M}$ on $(\mathbb{Z}/M\mathbb{Z})^{*}$.

The sequence $x_{t_{k,M}(x)},\cdots,x_{s_{k,M}(x)-1}$, ordered up to a cyclic shift, is called a \textit{cycle}.
The \textit{cycle length} is $c_{k,M}(x)=s_{k,M}(x)-t_{k,M}(x)$. The elements in the cycle are called \textit{periodic
 points} and their \textit{periods} are $c_{k,M}(x)$. So we can define a \textit{cycle length function}
 $c_{k,M}$ on $(\mathbb{Z}/M\mathbb{Z})^{*}$. In particular, \cite{FPS,KP} gave lower bounds for the largest period.

We denote by $P_{r}(k,M)$ and $P(k,M)$ respectively the number of periodic points with period $r$ and
the number of periodic points in $(\mathbb{Z}/M\mathbb{Z})^{*}$. Also, we denote by
$C_{r}(k,M)$ and $C(k,M)$ respectively the number of cycles with length $r$ and
the number of cycles in $(\mathbb{Z}/M\mathbb{Z})^{*}$. We denote the average values of $c_{k,M}(x)$ and $t_{k,M}(x)$ over all
$x\in (\mathbb{Z}/M\mathbb{Z})^{*}$ by $c(k,M)$ and $t(k,M)$ respectively,
\begin{equation}
c(k,M)=\frac{1}{\varphi(M)}\sum\limits_{x\in (\mathbb{Z}/M\mathbb{Z})^{*}}c_{k,M}(x),\qquad
t(k,M)=\frac{1}{\varphi(M)}\sum\limits_{x\in (\mathbb{Z}/M\mathbb{Z})^{*}}t_{k,M}(x),
\notag
\end{equation}
where $\varphi$ is the Euler totient function.

When $M$ is an odd prime power, we will derive explicit formulas for $P_{r}(k,M)$ and $C_{r}(k,M)$ by the results in
\cite{SK3}, and we will also derive explicit formulas for $c(k,M)$ and $t(k,M)$ which generalize those in \cite{VS}.

For two integers $r,m\ge 1$, we call the limit of $\lim\limits_{X \to \infty}\frac{1}{\pi(X)}\sum\limits_{p \le X}P_{r}(k,p^{m})$ the asymptotic mean number of periodic points with period $r$ in $(\mathbb{Z}/p^{m}\mathbb{Z})^{*}$ for different choices of prime $p$,
and we denote it by $AP_{r}(k,m)$. Similarly, we can define the asymptotic mean number for cycles with length $r$ and denote it by $AC_{r}(k,m)$.
We will derive explicit formulas for $AP_{r}(k,m)$ and $AC_{r}(k,m)$.

For an integer $m \ge 1$, following \cite{VS}, we study the average values of $P(k,p^{m})$ and $t(k,p^{m})$ over all
primes $p\le N$,
\begin{equation}
S_{0}(k,m,N)=\frac{1}{\pi(N)}\sum\limits_{p\le N}P(k,p^{m}),\qquad
S(k,m,N)=\frac{1}{\pi(N)}\sum\limits_{p\le N}t(k,p^{m}).
\notag
\end{equation}
where, as usual, $\pi(N)$ is the number of primes $p\le N$. Following the method in \cite{CS},
we will get asymptotic formulas for $S_{0}(k,m,N)$ and $S(k,m,N)$.

\section{Preparations}
For two integers $l$ and $n$, we denote their greatest common divisor by gcd$(l,n)$.
For a positive integer $n$, we denote by $\tau(n)$ the number of its positive divisors.
Theorem 4.9 in \cite{N} tells us that
\begin{equation}
\lim\limits_{X \to \infty}\frac{1}{\pi(X)}\sum\limits_{p \le X}{\rm gcd}(p-1,n)=\tau(n).
\end{equation}

For two integers $m\ge 1$ and $n\ge 2$, we denote the largest prime divisor of $n$ by $q$.
Then we have
\begin{equation}\label{p1}
\begin{array}{lll}
&&\lim\limits_{X \to \infty}\frac{1}{\pi(X)}\sum\limits_{p \le X}{\rm gcd}(p^{m-1}(p-1),n)\\
&=&\lim\limits_{X \to \infty}\frac{1}{\pi(X)}\sum\limits_{q<p \le X}{\rm gcd}(p^{m-1}(p-1),n)\\
&=&\lim\limits_{X \to \infty}\frac{1}{\pi(X)}\sum\limits_{q<p \le X}{\rm gcd}(p-1,n)\\
&=&\tau(n).
\end{array}
\end{equation}
Notice that if $p$ is an odd prime,
gcd$(p^{m}-p^{m-1},n)$ is the number of solutions of the equation $x^{n}=1$ in $(\mathbb{Z}/p^{m}\mathbb{Z})^{*}$.

Given two integers $a$ and $n$ with gcd$(a,n)=1$, following the method in the proof of Formula (2) in \cite{CS}, we can get
\begin{equation}
\sum\limits_{\substack{p\le X \\
p\equiv a\ ({\rm mod}\ n)}}p^{m}=\frac{X^{m+1}}{(m+1)\varphi(n){\rm ln}X}+O(X^{m+1}{\rm ln}^{-2}X).
\end{equation}
Then we have
\begin{equation}
\sum\limits_{\substack{p\le X \\
p\equiv a\ ({\rm mod}\ n)}}p^{m-1}(p-1)=\frac{X^{m+1}}{(m+1)\varphi(n){\rm ln}X}+O(X^{m+1}{\rm ln}^{-2}X).
\end{equation}

Following the same method in the proof of Formula (4) in \cite{CS}, we have
\begin{equation}
\sum\limits_{\substack{p\le X \\
p\equiv  a\ ({\rm mod}\ n)}}p^{m-1}(p-1)=O(\frac{X^{m+1}}{n}+X^{m}).
\end{equation}

\section{Main Results}                                                             

For two integers $d$ and $n$ satisfying gcd$(d,n)=1$, we denote the multiplicative order of $n$ modulo $d$ by ord$_{d}n$.
For an integer $n$ and a prime $p$, we denote $v_{p}(n)$ the exact power of $p$ dividing $n$.

Let $\mu$ be the M\"obius function. For a real number $a$, we denote $\lceil a \rceil$ the least integer which is not less than $a$.

Write $k=p_{1}^{n_{1}}p_{2}^{n_{2}}\cdots p_{s}^{n_{s}}\ge 2$,
where $p_{1},\cdots,p_{s}$ are distinct primes, $p_{1}<p_{2}<\cdots<p_{s}$ and
$n_{1},\cdots,n_{s}\ge 1$. Let $m$ be a fixed positive integer.

\begin{proposition}\label{pro1}
Let $p$ be an odd prime and $r$ be a positive integer.
Write $p^{m}-p^{m-1}=p_{1}^{r_{1}}\cdots p_{s}^{r_{s}}\cdot\rho$, where $r_{1},\cdots,r_{s}\ge 0$ are integers and
{\rm gcd}$(p_{1}\cdots p_{s},\rho)=1$. We have\\
{\rm (1)}
$C_{r}(k,p^{m})=\frac{1}{r}\sum\limits_{d|r}\mu(d){\rm gcd}(p^{m}-p^{m-1},k^{r/d}-1).$\\
{\rm (2)}
$P_{r}(k,p^{m})=\sum\limits_{d|r}\mu(d){\rm gcd}(p^{m}-p^{m-1},k^{r/d}-1).$\\
{\rm(3)} $P(k,p^{m})=\rho$.\\
{\rm(4)} $C(k,p^{m})=\sum\limits_{d|\rho}\frac{\varphi(d)}{{\rm ord}_{d}k}$.\\
{\rm(5)} For any $x\in (\mathbb{Z}/p^{m}\mathbb{Z})^{*}$, denote ${\rm ord}_{p^{m}}x$ by {\rm ord}x,
$c_{k,p^{m}}(x)={\rm ord}_{{\rm gcd}({\rm ord}x,\rho)}k$.\\
{\rm(6)} $c(k,p^{m})=\frac{1}{\rho}\sum\limits_{d|\rho}\varphi(d){\rm ord}_{d}k$.\\
{\rm(7)} For any $x\in (\mathbb{Z}/p^{m}\mathbb{Z})^{*}$, denote ${\rm ord}_{p^{m}}x$ by {\rm ord}x,
       \begin{equation}
              t_{k,p^{m}}(x)={\rm max}\left\{\lceil\frac{v_{p_{1}}({\rm ord}x)}{n_{1}}\rceil,
\lceil\frac{v_{p_{2}}({\rm ord}x)}{n_{2}}\rceil,\cdots,\lceil\frac{v_{p_{s}}({\rm ord}x)}{n_{s}}\rceil\right\}.
        \notag
        \end{equation}\\
{\rm(8)} $t(k,p^{m})=\frac{1}{p_{1}^{r_{1}}\cdots p_{s}^{r_{s}}}\sum\limits_{d|p_{1}^{r_{1}}\cdots p_{s}^{r_{s}}}\varphi(d){\rm max}\left\{\lceil\frac{v_{p_{1}}(d)}{n_{1}}\rceil,
\cdots,\lceil\frac{v_{p_{s}}(d)}{n_{s}}\rceil\right\}.$.
\end{proposition}
\begin{proof}
(1) and (2) By M\"obius inversion formula and Theorem 5.6 in \cite{SK3}.

(3) A special case of Corollary 3 in \cite{W}.

(4) By Theorem 2 and Theorem 3 in \cite{W}.

(5) By Lemma 3 and Theorem 2 in \cite{W}.

(6) Denote $p_{1}^{r_{1}}\cdots p_{s}^{r_{s}}$ by $w$, from (5), we have
\begin{equation}
\begin{array}{lll}
c(k,p^{m})&=&\frac{1}{p^{m}-p^{m-1}}\sum\limits_{x\in (\mathbb{Z}/p^{m}\mathbb{Z})^{*}}c_{k,p^{m}}(x)\\
&=& \frac{1}{p^{m}-p^{m-1}}\sum\limits_{d|\rho}\sum\limits_{n|w}\varphi(dn){\rm ord}_{d}k\\
&=&\frac{1}{p^{m}-p^{m-1}}\sum\limits_{n|w}\varphi(n)\sum\limits_{d|\rho}\varphi(d){\rm ord}_{d}k
=\frac{1}{\rho}\sum\limits_{d|\rho}\varphi(d){\rm ord}_{d}k.
\end{array}
\notag
\end{equation}

(7) Let $w_{x}$ be the factor of ord$x$ such that $\frac{{\rm ord}x}{w_{x}}$
is the largest factor relatively prime to $k$. By Lemma 3 in \cite{W}, we have
$t_{k,p^{m}}(x)$ is the least non-negative integer $l$ such that $w_{x}|k^{l}$.
In other words, $t_{k,p^{m}}(x)$ is the least non-negative integer $l$ such that
$v_{p_{i}}({\rm ord}x)\le ln_{i}$, for any $1\le i\le s$. Then we get the desired result.

(8) Notice that for any  $x\in (\mathbb{Z}/p^{m}\mathbb{Z})^{*}$, ord$x|(p^{m}-p^{m-1})$, and there are
$\varphi({\rm ord}x)$ elements with the order ord$x$. By (7), we have
\begin{equation}
t(k,p^{m})=\frac{1}{p^{m}-p^{m-1}}\sum\limits_{d|(p^{m}-p^{m-1})}\varphi(d){\rm max}\left\{\lceil\frac{v_{p_{1}}(d)}{n_{1}}\rceil,
\lceil\frac{v_{p_{2}}(d)}{n_{2}}\rceil,\cdots,\lceil\frac{v_{p_{s}}(d)}{n_{s}}\rceil\right\}.
\notag
\end{equation}
Furthermore, we have
\begin{equation}
\begin{array}{lll}
t(k,p^{m})&=&\frac{1}{p^{m}-p^{m-1}}\sum\limits_{d|p_{1}^{r_{1}}\cdots p_{s}^{r_{s}}\rho}\varphi(d){\rm max}\left\{\lceil\frac{v_{p_{1}}(d)}{n_{1}}\rceil,
\cdots,\lceil\frac{v_{p_{s}}(d)}{n_{s}}\rceil\right\}\\
&=&\frac{1}{p^{m}-p^{m-1}}\sum\limits_{i_{1}=0}^{r_{1}}\cdots \sum\limits_{i_{s}=0}^{r_{s}}\sum\limits_{d|\rho}
\varphi(p_{1}^{i_{1}}\cdots p_{s}^{i_{s}}d){\rm max}\left\{\lceil\frac{i_{1}}{n_{1}}\rceil,
\cdots,\lceil\frac{i_{s}}{n_{s}}\rceil\right\}\\
&=&\frac{1}{p^{m}-p^{m-1}}\sum\limits_{d|\rho}\varphi(d)\sum\limits_{i_{1}=0}^{r_{1}}\cdots \sum\limits_{i_{s}=0}^{r_{s}}
\varphi(p_{1}^{i_{1}}\cdots p_{s}^{i_{s}}){\rm max}\left\{\lceil\frac{i_{1}}{n_{1}}\rceil,
\cdots,\lceil\frac{i_{s}}{n_{s}}\rceil\right\}\\
&=&\frac{1}{p_{1}^{r_{1}}\cdots p_{s}^{r_{s}}}\sum\limits_{d|p_{1}^{r_{1}}\cdots p_{s}^{r_{s}}}\varphi(d){\rm max}\left\{\lceil\frac{v_{p_{1}}(d)}{n_{1}}\rceil,
\cdots,\lceil\frac{v_{p_{s}}(d)}{n_{s}}\rceil\right\}.\\
\end{array}
\notag
\end{equation}
\end{proof}

\begin{remark}
{\rm If we put $k=2$ and $m=1$, then the formulas (3),(4),(6) and (8) correspond to Theorem 6 in \cite{VS}.
}
\end{remark}

\begin{remark}
{\rm Since the conclusions in \cite{SK3} and \cite{W} we apply are about the general case of modulo a positive integer,
 it is easy to get similar formulas for the case of $p=2$.
}
\end{remark}

\begin{proposition}
Let $r$ be a positive integer, we have
\begin{equation}
 AP_{r}(k,m)=\sum\limits_{d|r}\mu(d)\tau(k^{r/d}-1),
\end{equation}
\begin{equation}
 AC_{r}(k,m)=\frac{1}{r}\sum\limits_{d|r}\mu(d)\tau(k^{r/d}-1).
\end{equation}
\end{proposition}
\begin{proof}
Combing (\ref{p1}) and Proposition \ref{pro1} (1) and (2), we can get the desired formulas.
\end{proof}

In the following, we denote by $\Omega$ the set of positive $\mathcal{S}$-units with $\mathcal{S}=\{p_{1},\cdots,p_{s}\}$.
Here a positive $\mathcal{S}$-unit means a positive integer whose prime divisors all belong to $\mathcal{S}$.

\begin{proposition}\label{f1}
We have
$\lim\limits_{N\to \infty}\frac{S_{0}(k,m,N)}{N^{m}}=\frac{1}{m+1}\left(\prod\limits_{i=1}^{s}\frac{p_{i}^{2}}{p_{i}^{2}-1}-1\right).$
\end{proposition}
\begin{proof}
Put $Q=p_{1}p_{2}\cdots p_{s}$ and denote by $\mathcal{U}_{Q}$ the set of integer $u,1\le u\le Q$, such that
gcd$(u,Q)=1$.

For each odd prime $p$, let $\rho_{p}$ be the largest divisor of $p^{m}-p^{m-1}$ coprime to $p_{1}p_{2}\cdots p_{s}$.
It is easy to see
\begin{equation}
\lim\limits_{N\to \infty}\frac{S_{0}(k,m,N)}{N^{m}}
=\lim\limits_{N\to \infty}\frac{1}{N^{m}\pi(N)}\sum\limits_{p_{s}<p\le N}\rho_{p}.
\notag
\end{equation}

Notice that if a prime $p>p_{s}$, then $v_{p_{i}}(p^{m}-p^{m-1})=v_{p_{i}}(p-1)$ for any $1\le i\le s$.
Hence, following the method in Theorem 2 of \cite{CS}, we have
\begin{equation}
\lim\limits_{N\to \infty}\frac{S_{0}(k,m,N)}{N^{m}}=\lim\limits_{N\to \infty}\frac{1}{N^{m}\pi(N)}\sum\limits_{q\in \Omega}q^{-1}\sum\limits_{u\in \mathcal{U}_{Q}}
\sum\limits_{\substack{p\le N \\
p\equiv qu+1\ ({\rm mod}\ qQ)}}(p^{m}-p^{m-1}).
\notag
\end{equation}

Following the method in Theorem 2 of \cite{CS}, we have
\begin{equation}
\lim\limits_{N\to \infty}\frac{S_{0}(k,m,N)}{N^{m}}=\frac{1}{m+1}\sum\limits_{q\in \Omega}\frac{1}{q^{2}}.
\notag
\end{equation}

Moreover, we have
\begin{equation}
\begin{array}{lll}
\sum\limits_{q\in \Omega}\frac{1}{q^{2}}&=&\sum\limits_{i_{1},\cdots,i_{s}=0}^{\infty}\frac{1}{(p_{1}^{i_{1}}\cdots p_{s}^{i_{s}})^{2}}-1\\
&=&\sum\limits_{i_{1}=0}^{\infty}\frac{1}{p_{1}^{2i_{1}}}\cdots \sum\limits_{i_{s}=0}^{\infty}\frac{1}{p_{s}^{2i_{s}}}-1\\
&=&\prod\limits_{i=1}^{s}\frac{p_{i}^{2}}{p_{i}^{2}-1}-1.
\end{array}
\notag
\end{equation}

Hence, we get the desired result.
\end{proof}

\begin{corollary}
We have $\frac{1}{k^{2}(m+1)}<\lim\limits_{N\to \infty}\frac{S_{0}(k,m,N)}{N^{m}}<\frac{2^{s}-1}{m+1}$.
\end{corollary}
\begin{proof}
Notice that for any prime $p$, we have 
$$
1+p^{-2}<\frac{p^{2}}{p^{2}-1}=1+\frac{1}{p^{2}-1}<2.
$$
\end{proof}

Given $q=p_{1}^{r_{1}}\cdots p_{s}^{r_{s}}\in \Omega$, we denote
\begin{equation}
\psi(q)=\frac{1}{q}\sum\limits_{d|q}\varphi(d){\rm max}\left\{\lceil\frac{v_{p_{1}}(d)}{n_{1}}\rceil,
\cdots,\lceil\frac{v_{p_{s}}(d)}{n_{s}}\rceil\right\}.
\notag
\end{equation}

\begin{proposition}
We have $\lim\limits_{N\to \infty}S(k,m,N)=\sum\limits_{q\in \Omega}\frac{\psi(q)}{q}.$
\end{proposition}
\begin{proof}
Given $q=p_{1}^{r_{1}}\cdots p_{s}^{r_{s}}\in \Omega$. Suppose $r_{1}\ge 1$, we want to estimate $\frac{1}{q}\sum\limits_{d|q}\varphi(d)\lceil\frac{v_{p_{1}}(d)}{n_{1}}\rceil$.
For simplicity, we replace $p_{1},r_{1}$ and $n_{1}$ by $p,r$ and $n$ respectively.
By division algorithm, we write $r=ln+d$ with $0\le d<n$. We have
\begin{equation}
\begin{array}{lll}
\frac{1}{q}\sum\limits_{d|q}\varphi(d)\lceil\frac{v_{p}(d)}{n}\rceil
&=&\frac{1}{p^{r}}\sum\limits_{d|p^{r}}\varphi(d)\lceil\frac{v_{p}(d)}{n}\rceil\\
&=&\frac{p-1}{p^{r}}\sum\limits_{i=1}^{r}p^{i-1}\lceil\frac{i}{n}\rceil\\
&=&\frac{p-1}{p^{r}}\left[\sum\limits_{i=1}^{n}p^{i-1}+\sum\limits_{i=n+1}^{2n}2p^{i-1}
+\cdots+\sum\limits_{i=(l-1)n+1}^{ln}lp^{i-1}+\sum\limits_{i=ln+1}^{ln+d}(l+1)p^{i-1}\right]\\
&=&\frac{p^{n}-1}{p^{r}}\left[1+2p^{n}+\cdots+lp^{(l-1)n}\right]+\frac{(l+1)p^{ln}(p^{d}-1)}{p^{r}}\\
&=&\frac{lp^{ln}}{p^{r}}-\frac{p^{ln}-1}{p^{r}(p^{n}-1)}+\frac{(l+1)p^{ln}(p^{d}-1)}{p^{r}}\\
&\le&l+(l+1)\le 3r.
\end{array}
\notag
\end{equation}
Hence, we have
\begin{equation}\label{bound}
\begin{array}{lll}
\psi(q)&\le& \frac{1}{q}\sum\limits_{d|q}\varphi(d)\left(\lceil\frac{v_{p_{1}}(d)}{n_{1}}\rceil+
\cdots+\lceil\frac{v_{p_{s}}(d)}{n_{s}}\rceil\right)\\
&\le&3(r_{1}+\cdots+r_{s})\\
&\le&\frac{3}{\ln 2}\ln q=O({\rm ln}q).\\
\end{array}
\end{equation}

Similarly to Proposition \ref{f1}, by Proposition \ref{pro1} (8), we have
\begin{equation}
\lim\limits_{N\to \infty}S(k,m,N)=\lim\limits_{N\to \infty}\frac{1}{\pi(N)}\sum\limits_{q\in \Omega}\psi(q)\sum\limits_{u\in \mathcal{U}_{Q}}
\sum\limits_{\substack{p\le N \\
p\equiv qu+1\ ({\rm mod}\ qQ)}}1.
\notag
\end{equation}
Then following the method in Theorem 2 of \cite{CS}, we can get the desired result.
\end{proof}

\begin{corollary}
We have $\frac{1}{k}<\lim\limits_{N\to \infty}S(k,m,N)<\frac{5\sqrt{p_{1}}\cdots\sqrt{p_{s}}}{(\sqrt{p_{1}}-1)\cdots (\sqrt{p_{s}}-1)}$.
\end{corollary}
\begin{proof}
On one hand we have
\begin{equation}
\begin{array}{lll}
\sum\limits_{q\in \Omega}\frac{\psi(q)}{q}&>&\sum\limits_{i_{1}\ge n_{1},\cdots,i_{s}\ge n_{s}}\frac{\varphi(p_{1}^{i_{1}}\cdots p_{s}^{i_{s}})}{(p_{1}^{i_{1}}\cdots p_{s}^{i_{s}})^{2}}\\
&=&\frac{(p_{1}-1)\cdots (p_{s}-1)}{p_{1}\cdots p_{s}}\sum\limits_{i_{1}\ge n_{1}}^{\infty}\frac{1}{p_{1}^{i_{1}}}\cdots \sum\limits_{i_{s}\ge n_{s}}^{\infty}\frac{1}{p_{s}^{i_{s}}}\\
&=&\frac{1}{k}.
\end{array}
\notag
\end{equation}
On the other hand, by (\ref{bound}) we have $\psi(q)<5\ln q$, then we have
\begin{equation}
\begin{array}{lll}
\sum\limits_{q\in \Omega}\frac{\psi(q)}{q}&<&\sum\limits_{q\in \Omega}\frac{5\ln q}{q}\\
&<& 5\sum\limits_{q\in \Omega}\frac{1}{\sqrt{q}}\\
&=&5\sum\limits_{i_{1}=0,\cdots,i_{s}=0}\frac{1}{\sqrt{p_{1}^{i_{1}}\cdots p_{s}^{i_{s}}}}\\
&=&\frac{5\sqrt{p_{1}}\cdots\sqrt{p_{s}}}{(\sqrt{p_{1}}-1)\cdots (\sqrt{p_{s}}-1)}.\\
\end{array}
\notag
\end{equation}

\end{proof}

\section{Remarks on the general case}
In this section, we will give some remarks on the case of modulo a positive integer.

We can deduce formulas for $C_{r}(k,M)$ and $P_{r}(k,M)$ directly from
Theorem 5.6 in \cite{SK3} . Corollary 3 in \cite{W} has given a formula for $P(k,M)$.
We can also derive a formula for $C(K,M)$ directly by applying Theorem 2 and Theorem 3 in \cite{W}.

Following the same methods, we can easily determine the cycle length function $c_{k,M}(x)$ and the tail function $t_{k,M}(x)$ on $(\mathbb{Z}/M\mathbb{Z})^{*}$,
then we can get formulas for $c(k,M)$ and $t(k,M)$.

In fact, \cite{W} and \cite{SK3} can tell us more information about the properties of repeated
exponentiation modulo a positive integer.

\section{Acknowledgment}
We would like to thank Prof. I.E. Shparlinski for suggesting
this problem and for his helpful advice.

\bibliographystyle{amsplain}

\begin{thebibliography}{20}

\bibitem{BBS} L. Blum, M. Blum and M. Shub, \textit{A simple unpredictable pseudo-random
number generator}, SIAM J. Comp., \textbf{15} (1986), 364-383.

\bibitem{CS} W.-S. Chou and I.E. Shparlinski, \textit{On the Cycle Structure of Repeated Exponentiation Modulo a Prime}, J. Number Theory,
\textbf{107} (2004), 345-356.

\bibitem{E} E. D. El-Mahassni, \textit{On the Distribution of the Power Generator over a Residue Ring for Parts of the Period}, Rev. Mat. Complut.,
 \textbf{21} (2008), 319-325.

\bibitem{FPS} J. B. Friedlander, C. Pomerance and I. E. Shparlinski, \textit{Period of the
power generator and small values of Carmichael's
function}, Math. Comp.,  \textbf{70} (2001), 1591-1605.

\bibitem{KP} P. Kurlberg and C. Pomerance, \textit{On the period of the linear congruential
and power generators}, Acta Arith., \textbf{119} (2005), 149-169.

\bibitem{L} J. C. Lagarias, \textit{Pseudorandom number generators in cryptography and
number theory}, Proc. Symp. in Appl. Math., Amer. Math. Soc., Providence,
RI, \textbf{42} (1990), 115-143.

\bibitem{N} M. Nilsson, \textit{ Cycles of monomials and perturbated monomial p-adic
dynamical systems}, Ann. Math. Blaise Pascal, {\bf 7}(1) (2000), 37-63.

\bibitem{SK1}L. Somer, M. K$\check{\rm r}\acute{\i}\check{\rm z}$ek, \textit{Structure of digraphs associated with quadratic congruences with composite moduli}, Discrete Math. \textbf{306} (2006) 2174-2185.

\bibitem{SK2}L. Somer, M. K$\check{\rm r}\acute{\i}\check{\rm z}$ek, \textit{On semiregular digraphs of the congruence $x^{k}\equiv y$ (mod n)}, Comment. Math. Univ. Carolin. \textbf{48} (2007), 41-58.

\bibitem{SK3}L. Somer, M. K$\check{\rm r}\acute{\i}\check{\rm z}$ek, \textit{On symmetric digraphs of the congruence $x^{k}\equiv y$ (mod n)},
  Discrete Math. \textbf{309} (2009) 1999-2009.

\bibitem{VS}T. Vasiga, J. Shallit, \textit{On the iteration of certain quadratic maps over $GF(p)$}, Discrete Math. \textbf{277} (2004) 219-240.

\bibitem{W}B. Wilson, \textit{Power digraphs modulo n}, Fibonacci Quart. \textbf{36} (1998) 229-239.

\end{thebibliography}

\end{document}